\documentclass[10pt]{amsart}
\usepackage{amsmath,amscd,latexsym,verbatim,amssymb}
\usepackage{times}       
 \newcommand{\resp}{{\it resp.} }
\newcommand{\cf}{{\it cf.} }
\newcommand{\ie}{{\it i.e.} }
\newcommand{\eg}{{\it e.g.} }

\newcommand{\loccit}{{\it loc. cit.} }

\newcommand{\R}{\mathbf{R}}
\newcommand{\C}{\mathbf{C}}

  \newcommand{\sA}{{\mathcal{A}}}

\newcommand{\inj}{\hookrightarrow}

\newcommand{\surj}{\rightarrow\!\!\!\!\!\rightarrow}

\renewcommand{\epsilon}{\varepsilon}
\renewcommand{\phi}{\varphi}

\renewcommand{\lim}{\varprojlim}

\font\smit= cmti10 at 10pt

\newcounter{spec}
{\end{list}}

\swapnumbers

\newtheorem{thm}{Theorem}[subsection]

\theoremstyle{definition}

\newtheorem{rem}[thm]{Remark}
\newtheorem{rems}[thm]{Remarks}

\numberwithin{equation}{section}

 at10pt
 at12pt

\renewcommand{\qed}{\hfill $\square$\medskip}

\begin{document}

\title{On relative integral monodromy of abelian logarithms and normal functions.} 
\author{Yves
Andr\'e}
  \address{Institut Math\'ematique de Jussieu, Sorbonne Universit\'e,
  Paris 05\\France.}
\email{yves.andre@imj-prg.fr}
  \date{8/10/2024}

\maketitle
 
  \begin{abstract}  The relative algebraic monodromy of abelian logarithms (defined as the kernel of a map between algebraic monodromy groups attached to an abelian scheme with and without a section) was computed in \cite{A1}: under natural assumptions, this vector group turns out to be maximal. 
  
  The relative integral monodromy of abelian logarithms is defined similarly as a kernel of integral monodromy groups, without taking Zariski closures. We show that if the integral monodromy of the abelian scheme is a lattice in the algebraic monodromy (which is not always the case), then the relative integral monodromy of the abelian logarithm is also a lattice in the relative algebraic monodromy. The proof uses a Hodge-theoretic interpretation of sections of abelian schemes.
  
  We also consider relative integral monodromy groups in the more general context of normal functions.   
  \end{abstract}



\begin{sloppypar} 
  
 \section{Introduction}  
\subsection{} The relative monodromy of abelian logarithms of the title refers to the following situation: let $f: A \to S$ be an abelian scheme parametrized by a connected smooth complex variety endowed with a section $\sigma$. The Betti realization of the family of $1$-motives $[\mathbb Z_S \stackrel{1\mapsto \sigma}\to A]$ lies in an exact sequence of local systems 
\begin{equation}\label{eq1} 0 \to V_{\mathbb Z}= R_1f^{ an}_\ast \mathbb Z  \to {\bf V}_{\mathbb Z} \to \mathbb Z_S\to 0 \end{equation}
 on $S^{an}$, where $\mathbb Z_S$ denotes the constant local system $\mathbb Z$ on $S^{an}$ \cf \eg \cite[\S 10]{D}\cite[\S 1.4]{A2}. 

Let $s\in S(\mathbb C)$ be a base point, $\Gamma$ (\resp $\bf \Gamma$) be the image of $\pi_1(S(\C),s)$ in $GL(V_s)$ (\resp $GL({\bf V}_s)$), and $H$ (\resp $\bf H$) be its Zariski-closure. By a theorem of Deligne, the neutral component of $H$ is a semisimple group \cite{D1}. 

Let $\Delta$ (the {\it relative integral monodromy}) be the kernel of ${\bf \Gamma} \to \Gamma$ and $U$ (the {\it relative algebraic monodromy}) be the kernel of ${\bf H} \to H$, so that $U$ is a vector group and $\Delta$ is a subgroup of $U(\mathbb Z)$.  
 In \cite{A1} (\cf also \cite[\S\S 1.4, 1.5]{A2}), the following result was proven:

\begin{thm} Assume that 
\item $i)$ $A/S$ has no fixed part, even after pull-back to a finite etale cover of $S$,
\item $ii)$ $\mathbb Z.\sigma(S) $ is Zariski-dense in $A$ (\ie the image of the section is not contained in any proper $S$-subgroup of $A$). 

\smallskip Then there is a natural isomorphism $U(\mathbb Z)\cong H_1(A_s, \mathbb Z) $.  
\end{thm}

This has been used in a number of contexts, notably in the study of Betti maps in Diophantine Geometry \cite{ACZ}\footnote{the name ``Betti map", due to D. Bertrand, actually refers to this context of Betti realizations.}.

\subsection{} It is tempting to replace the conclusion of Theorem 1.1.1 by\footnote{\cf \eg \cite[\S 3.1]{V} - which however only uses 1.1.1, not its overinterpretation in terms of integral relative monodromy.}: 

\smallskip   {\it $\;\;\Delta$ is of finite index in $H_1(A_s, \mathbb Z) \; \;\; (?) $  } 

\smallskip However, as P. Corvaja and U. Zannier pointed out \cite{CZ}, it is not clear whether $\Delta$ is Zariski-dense in $U$: indeed, taking a kernel does not commute with taking a Zariski-closure.
 
  The aim of this note is to clarify this point. We prove: 

\begin{thm} Assume moreover that 

\item $iii)$ $\Gamma$ is a lattice in $H_{\mathbb R}$. 

\smallskip Then $\Delta$ is a lattice in $U(\mathbb Z) \,$ (equivalently, $\bf\Gamma$ is a lattice in ${\bf H}_{\mathbb R}$).  
\end{thm}

Assumption $iii)$ is equivalent to: $\Gamma $ is of finite index in $H(\mathbb Z)$ \cite{BH}. It holds for instance in the situation of elliptic schemes, which is settled in \cite{CZ}\footnote{in this particular situation, Theorem 1.1.1 was obtained previously by D. Bertrand, and we reformulate \cite{CZ} in Hodge-theoretic terms in \S 2 below. Theorem 1.2.1 in the case of polarized abelian schemes of relative dimension $g$ over a base which is finite surjective over the moduli space $\mathcal A_g$ was considered in \cite{DT}).}. 


 However, $iii)$ does not always hold: there are abelian schemes with non-lattice monodromy $\Gamma$, \cf \cite{N}. We do not know whether the conclusion holds in such cases. 

\subsection{} One key point in our proof is the Hodge-theoretic interpretation of (holomorphic versus algebraic) sections of abelian schemes, which is of independent interest (\cf Theorem 3.2.1, which depends on a theorem of Graber-Starr). 

This allows us to reduce the question to the ``modular case" and apply the theorem of Mok-To on sections of Kuga fibrations.

\subsection{} In most cases, the study of relative integral monodromy can be settled by general results on vanishing of the first cohomology group $H^1(\Gamma, V)$ of representations of integral monodromy groups, when the latter are arithmetic. We review these results (Margulis-Starkov, Bass-Milnor-Serre) in \S\S 5.1, 5.2. 
 Such arguments, which do not bear on the geometry of the situation, allow to deal with relative integral monodromy in the more general situation of normal functions, notably in the context of Lefschetz pencils \S 5.3. 

When these general results do not apply ($\Gamma$ is arithmetic but $H^1(\Gamma, V)\neq 0$ $-$ as in the case of elliptic families where $\Gamma$ is an arithmetic subgroup of $SL(2)$), our proof rests on the fact that $H^1(\Gamma, V)$ carries a mixed Hodge structure with $H^1(\Gamma, V)^{0,0}=0.$

\section{The case of families of elliptic curves: Hodge-theoretic interpretation of \cite{CZ}}

\subsection{} Replacing $S$ by a finite etale cover, we may assume that for some $N>2$, the $N$-torsion of $f$ is defined over $S$, whence a (non-constant) map $S \to \mathcal A_{1,N} = X(N)$ from $S$ to the modular curve, $\Gamma$ being a subgroup of finite index of $\Gamma(N)$. Denoting by $S'$ the finite etale cover of $X(N)$ with fundamental group $\Gamma$, there is a cartesian square    $$\begin{matrix}    A  &  \stackrel{h}\to  & A'\\   f \downarrow&&  \downarrow {f'} \\ S &  \stackrel{g}\to & S' . \end{matrix}$$  

The cocycle $  e\in H^1(S^{an}, R_1f^{an}_\ast \mathbb Q )= H^1(\pi_1(S, s), H_1(A_s, \mathbb Q ))$  attached to the section $\sigma$ comes (via $(h\sigma)^\ast$) from a cocycle in $H^1({\bf \Gamma}, H_1(A_s, \mathbb Q )) $.
 
\smallskip  The rank of $\Delta$ is a priori $0,1$ or $2$.
 If it is $2$, $\Delta$ is a lattice in $ H_1(A_s, \mathbb R) $. 
 
  It cannot be $1$, because $\mathbb Q. \Delta$ is monodromy-invariant, and $\Gamma$-invariant subspaces of $ H_1(A_s, \mathbb Q) $ are $0$ or $ H_1(A_s, \mathbb Q) $.   
 
 If it is $0$, \ie  ${\bf \Gamma} = \Gamma$, then $e$ comes via $g^\ast = (f'h\sigma)^\ast$ from a cocycle in $$ H^1(\Gamma,  H_1(A_s, \mathbb Q )) = H^1(\pi_1(S', g(s)), H_1(A'_{h(s)}, \mathbb Q ))= H^1(S'^{an} , R_1{f}^{'an}_\ast \mathbb Q ).$$  
  
 The morphism 
 $$g^\ast: H^1(S'^{an} , R_1{ f}^{'an}_\ast \mathbb Q ) \to H^1(S^{an}, R_1f^{an}_\ast \mathbb Q )$$
  is a morphism of mixed Hodge structures. The morphism $g$ extends to a morphism $\bar g: \bar S \to \overline{S'}$ between compactifications $j: S \inj \bar S, {j'}: S' \inj \overline{S'}$, and the weight-$0$ part of $g^\ast$ can be identified with
 $$H^1(\overline{S'} , j'_\ast R^1{f'}^{an}_\ast \mathbb Q(1)) \to H^1(\bar S, j_\ast R^1f^{an}_\ast \mathbb Q(1)).$$
 According to \cite{Z}, the $(0,0)$-component of this map can be identified with the map between sections $A'(S')\otimes \mathbb Q \to A(S)\otimes \mathbb Q$ induced by $g$. 
  Therefore, if $\Delta =0$, then $\sigma$ comes from a section of $f'$. But $A'(S') $ is torsion according to Shioda's theorem \cite{Sh}. Hence $\sigma $ is torsion and $U=0$ in that case. 
  
  \begin{rem} In the case of Legendre pencil $f: A \to  S= S' = \mathbb P^1 \setminus \{0, 1, \infty \}$, 
  an example of a non-torsion class $e\in H^1(S, V_s)$ coming from a holomorphic, non-algebraic, section of $f^{an}$ is given in \cite[(2.6)]{CZ}. \end{rem}
 
\section{Hodge-theoretic interpretation of sections of abelian schemes}

Let $f: A\to S$ be an abelian scheme {\it without fixed part}. Let $V_\mathbb Z $ be the local system $ R_1f^{an}_\ast \mathbb Z$ on $S^{an}$ - its fiber $H_1(A^{an}, \mathbb Z)$ at $s\in S(\mathbb C)$ is a representation of $\Pi := \pi_1(S(\mathbb C), s)$), and it underlies a polarized integral variation of Hodge structures of type $(0,-1),(-1, 0)$. 

Let us discuss holomorphic sections versus algebraic sections of $A/S$ from a Hodge-theoretic viewpoint.

 \subsection{Holomorphic sections} Let ${Sec}\, f^{an}$ be the sheaf of germs of holomorphic sections of $f $. One has  $H^0(S^{an}, {Sec}\, f^{an})= \mathcal O_{A^{an}}(S^{an})$ and an exact sequence of sheaves on $S^{an}$:  
 \begin{equation}\label{eq2} 0\to V_{\mathbb Z} \to { Lie} \, A^{an} \to  {Sec}\, f^{an}\to 0.\end{equation}
 In the Hodge-theoretic viewpoint, it is useful to identify ${ Lie} \, A^{an}$  with the notch $ Gr^{-1}{\mathcal H}$ of the holomorphic vector bundle $\,\mathcal H :=  V\otimes \mathcal O_{S^{an}} $ with respect to the Hodge filtration.  
 Since $ H^0(S^{an}, V_{\mathbb Z} )=0$ by assumption (no fixed part),  one gets an exact sequence  
\begin{equation}\label{eq3} 0 \to  H^0(S^{an}, Gr^{-1}{\mathcal H}) \to H^0(S^{an},\mathcal \underline{Sec}\, f^{an}) \to H^1(S^{an}, V_{\mathbb Z}) \to H^1(S^{an}, Gr^{-1}{\mathcal H} ) .\end{equation}
 Elements of $H^1(S^{an}, V_{\mathbb Z})= H^1(\Pi, V_{\mathbb Z , s})$ parametrize extensions of the local systems $\mathbb Z$ by $V_{\mathbb Z}$: 
\begin{equation}\label{eq4} H^1(S^{an}, V_{\mathbb Z}) \cong {\rm{Ext}}^1_{loc. sys}(\mathbb Z_S, V_{\mathbb Z} ),\end{equation}
 while elements of $H^0(S^{an},\mathcal \underline{Sec}\, f^{an})$ (holomorphic sections of $A/S$) parametrize extensions of the variation of mixed Hodge structures of $\mathbb Z$ by $V_{\mathbb Z}$:
\begin{equation}\label{eq5} H^0(S^{an}, \mathcal \underline{Sec}\, f^{an}) \cong {\rm{Ext}}^1_{VMHS}(\mathbb Z_S, V_{\mathbb Z} )\end{equation}  
 (indeed: fiberwise, extensions of $\mathbb Z$ by $V_{\mathbb Z , s}$ as mixed Hodge structures are parametrized by points of $A_s$ \cite{Ca}; in the relative situation, note that Griffiths' transversality, which is part of the axioms of a VMHS, is automatic here\footnote{in the language of normal functions, any holomorphic section (which is the same as a normal function in this situation) is horizontal.} since there is only one intermediate step in the Hodge filtration of $\mathcal H$, namely $F^{0}$).
   
   The kernel of the map which forgets Hodge structures 
\begin{equation}\label{eq6}{\rm{Ext}}^1_{VMHS}(\mathbb Z_S, V_{\mathbb Z} ) \to {\rm{Ext}}^1_{loc. sys}(\mathbb Z_S, V_{\mathbb Z} )=  H^1(S^{an}, V_{\mathbb Z}) \end{equation}
 is therefore $H^0(S^{an}, Gr^{-1}{\mathcal H})$, via \eqref{eq3}, \eqref{eq5} and \eqref{eq4} (elements of $H^0(S^{an}, Gr^{-1}{\mathcal H})$ just perturb the Hodge filtration). 
   
 \subsection{Algebraic sections} Steenbrink-Zucker and Kashiwara have defined a notion of {\it admissible} variation of mixed Hodge structures ${\bf V}_{\mathbb Z}$ (with respect to a smooth compactification $\bar S$ of $S$ \cite{SZ}\cite{K} - a notion which can be tested by restriction to curves in $S$, and which is equivalent to the extendability of the VMHS  to a mixed Hodge module on $\bar S$ in M. Saito's sense).  

On the other hand, Deligne-Zucker-Cattani-Kaplan-Kashiwara showed that $H^1(S^{an}, V_{\mathbb Z}) $ underlies a canonical polarized mixed Hodge structure \cite{CK}\cite{K}. The restriction of \eqref{eq6} to admissible VMHS gives rise to an injection 
 \begin{equation}\label{eq7}   {\rm{Ext}}^1_{VMHS^{ad}}(\mathbb Z_S, V_{\mathbb Z}) \inj H^1(S^{an}, V_{\mathbb Z})^{0,0}    \end{equation}  
 (injectivity follows from the normality theorem of \cite{A1} for admissible VMHS, which itself relies on the theorem of the fixed part of Steenbrink-Zucker and Kashiwara, \cf \cite[\S 1.4]{A2}).
 
  \smallskip Any {\it algebraic} section $\sigma\in \mathcal O_A(S)$ gives rise to an integral admissible variation of mixed Hodge structures ${\bf V}_\mathbb Z \in  {\rm{Ext}}^1_{VMHS^{ad}}(\mathbb Z_S, V_{\mathbb Z})$. In particular, the class in $H^1(S^{an}, V_{\mathbb Z})$ of any algebraic section belongs to $H^1(S^{an}, V_{\mathbb Z})^{0,0}$, and one gets injections 
  \begin{equation}\label{eq8}  \mathcal O_A(S)  \inj {\rm{Ext}}^1_{VMHS^{ad}}(\mathbb Z_S, V_{\mathbb Z}  ) \inj H^1(S^{an}, V_\mathbb Z )^{0,0}\end{equation}
(injectivity of the left map can be checked fiberwise, where it reduces to the fact that the MHS associated to the $1$-motive $[\mathbb Z \stackrel{1\mapsto \sigma(s)}\to A_s]$ splits only if $\sigma(s)=0$ \cite[\S 10]{D}).    
 
\smallskip All the above discussion may be seen as a well-known special instance ``sections of abelian schemes" of the modern theory of normal functions initiated by Green and Griffiths.
 However, the following comparison between the three terms which intervene in \eqref{eq8} does not belong to the general theory but only to this special instance: indeed, it essentially relies upon a geometric theorem by Graber and Starr \cite[Theorem 1.3]{GS}. 
 
 We continue to assume that $A/S$ has no fixed part.

\begin{thm}\label{3.2.1} \eqref{eq8} induces isomorphisms
  \begin{equation}\label{eq9}  \mathcal O_A(S)\otimes \mathbb Q  \cong {\rm{Ext}}^1_{VMHS^{ad}}(\mathbb Q_S, V  ) \cong H^1(S^{an}, V )^{0,0}.\end{equation}
 Moreover, if either 
\smallskip  \item a) $S$ is proper or contained in a proper variety with boundary of codimension $>1$, or 
  \item  b) $S$ has a compactification with boundary with normal crossings and the local monodromies of $V$ at the boundary are unipotent, or
  \item c) if for some $N\geq 3$, the pull-back $A'/S'$ of $A/S$ over some Galois covering $S'/S$ where $A'$ acquires $N$-torsion has no fixed part, 
  
\smallskip  then the maps in \eqref{eq8} are isomorphisms. \end{thm} 

(Only \eqref{eq9} will be used in the sequel).
      
  \begin{proof}  Let $j: S\inj \bar S$ be an open immersion into a proper variety. Let $\pi: \bar S \to \mathbb P^d$ be a surjective generically finite etale morphism. If $d = \dim S >1$, the theorem of Graber-Starr says that for any general conic $C' \subset \mathbb P^d$, setting $C:= C'\times_{\mathbb P^d} S$ (which is a smooth curve $ C \stackrel{i}{\inj} S$), the restriction map on sections is an isomorphism
    \begin{equation}\label{eq10} \mathcal O_A(S) \cong \mathcal O_{A_C}(C). \end{equation}
    If $d=1$, we take $C=S$ in the sequel of the proof.
    Let us consider the commutative diagram
         \begin{equation}\label{eq11}\begin{matrix}  \mathcal O_A(S)  &\inj {\rm{Ext}}^1_{VMHS^{ad}}(\mathbb Z_S, V_{\mathbb Z}  ) &\inj H^1(S^{an}, V_\mathbb Z )^{0,0}   \\
            i^\ast  \downarrow&  i^\ast \downarrow &  \downarrow i^\ast  \\
             \mathcal O_{A_C}(C) & \inj {\rm{Ext}}^1_{VMHS^{ad}}(\mathbb Z_C, i^\ast V_{\mathbb Z}  ) &\inj H^1(C^{an}, i^\ast V_\mathbb Z )^{0,0}.  \end{matrix}\end{equation}
      The first vertical map \eqref{eq11} is an isomorphism. The third vertical map is injective since $\pi_1(C)\to \pi_1(S)$ is surjective by Bertini-Lefschetz, hence 
       \begin{equation}\label{eq12}H^1(S^{an}, V_\mathbb Z )\subset  H^1(C^{an}, i^\ast V_\mathbb Z ). \end{equation}   
      
      In case a), $C$ is proper, hence $\mathcal O_{A_C}(C) = \mathcal O_{A^{an}_C}(C^{an})$ and the bottom line of \eqref{eq11} is a sequence of isomorphims. And so is the top line.
      
      In case b), the bottom line of \eqref{eq11} is a also sequence of isomorphims: by \cite[10.2]{Z} (which uses the assumption: no fixed part), the map 
      $   \mathcal O_{A_C}(C) \to H^1(C^{an}, i^\ast V_\mathbb Z )^{0,0}$  restricts to an isomorphism 
     \begin{equation}\label{eq12}  \mathcal O^0_{A_C}(C) \cong H^1(\bar C^{an}, j_{C \ast} i^\ast V_\mathbb Z )^{0,0},  \end{equation} 
      where $\mathcal O^0_{A_C}(C)$ denotes the subgroup of sections which factor through the identity component of the N\'eron model, which is the case of all sections if the local monodromies are unipotent. Therefore the top line of \eqref{eq11} is again a sequence of isomorphisms.
      
      In case c), let $G$ be the Galois group of $S'/S$, and $\Pi'$ be the kernel of $\Pi\to G$. One has   $\mathcal O_A(S) \cong   \mathcal O_{A_{S'}}(S')^G$. Let us consider the commutative diagram (with exact bottom row)
         \begin{equation}\label{eq13}\begin{matrix}  &    \mathcal O_A(S) & \to   \mathcal O_{A_{S'}}(S')^G  \\
            &    \downarrow &  \downarrow &  \\  H^1(G, V_{\mathbb Z, s}^{\Pi'})  &\to   H^1(\Pi,  V_{\mathbb Z, s})  &\to H^1(\Pi',  V_{\mathbb Z, s})^G  \end{matrix}\end{equation} Since $V_{\mathbb Z, s}^{\Pi'}= 0$ by assumption (no fixed part), one has $H^1(G, V_{\mathbb Z, s}^{\Pi'})=0$, so that the map  $ H^1(\Pi,  V_{\mathbb Z, s})   \to H^1(\Pi',  V_{\mathbb Z, s})^G$ is injective.
  On the other hand, since $N\geq 3$, $A'/S'$ falls in case $b)$, so that $\mathcal O_{A_{S'}}(S') \cong H^1(\Pi',  V_{\mathbb Z, s})^{0,0}$, and one concludes that  $\mathcal O_A(S)\cong H^1(\Pi,  V_{\mathbb Z, s})^{0,0}$.  
      
   \smallskip     In general 
      \begin{equation}\label{eq14}   H^1(\bar C^{an}, j_{C \ast} i^\ast V )^{0,0}= H^1(C^{an}, i^\ast V )^{0,0}.   \end{equation} 
  Taking into account that 
  $  {\rm{Ext}}^1_{VMHS^{ad}}(\mathbb Z_S,  V_{\mathbb Z}  )\otimes \mathbb Q \cong {\rm{Ext}}^1_{VMHS^{ad}}(\mathbb Q_S, V  )$ \cite[Lemma 2.8]{BFZP}, one concludes \eqref{eq9} from \eqref{eq13}, \eqref{eq12}, \eqref{eq10}. \end{proof}
    
\begin{rems} 1) Even when $f$ extends to a semiabelian scheme $\tilde f$ over a compactification of $S$ with boundary with normal crossings and $\dim S >1$, neither $\sigma$ nor any non-zero multiple of $\sigma$ will extend to a section of $\tilde f$ in general: it may be necessary to blow up some subscheme of $S$ to achieve this (roughly, the condition of admissibility on a holomorphic section has to do with moderate growth, and further with meromorphy, not holomorphy, at infinity). Therefore, the inquiry about algebraicity of sections cannot be performed in general by looking at limits in various extensions of $A^{an}/S^{an}$ as in the general theory of normal functions: instead, we used the Graber-Starr theorem.

\smallskip \noindent 2) An element of ${\mathcal O}_A(S)\otimes \mathbb Q$ gives rise to a $1$-motive over $S$ up to isogeny.
Its rational Betti realization $V$ is an admissible variation of mixed Hodge structures \cite{A1}. Theorem 3.2.1 suggests the question of the converse: {\it does any polarizable admissible variation of mixed Hodge structures of types $(0,0), (0,-1), (-1, 0), (-1,-1)$ comes from a $1$-motive up to isogeny over $S$?} 

In this respect, note that by Cartier duality, Theorem 3.2.1 also provides a Hodge-theoretic interpretation of extensions of abelian schemes by $\mathbb G_m$. To settle the question, it remains to match these two interpretations via biextensions...

\smallskip \noindent 3) After finishing this note, the author received the preprint \cite{DT2} where the authors show, under mild technical assumptions, that for a {\it ramified} section $\sigma$, $\Delta$ is a lattice in $U(\mathbb R)$. This is especially interesting since it is not assumed that $\Gamma$ is a lattice in $H(\mathbb R)$. 

Assume for simplicity that $S$ is a curve. Assume that the torsion group $A[N]$ is constant for some $N\geq 3$, whence a map $S\stackrel{\phi}{\to} \sA_{g, N}$, and a finite extension $\mathbb C(S)/\mathbb C(\phi(S))$ (which we may assume to be Galois). Let $S'$ be the normalization of $\phi(S)$ in the maximal subextension $\mathbb C(S')\subset \mathbb C(S)$ unramified at all points of $\phi(S)$. There is a factorization $S\to S'\to \phi(S)$ and 
 $f$ comes from an abelian scheme $f': A'\to S'$.
  A section of $f$ is {\it ramified} if it does not come from a section of $f'$.   

In order to prove that $\Delta\neq 0$ for ramified sections, \cite{DT2} uses the Betti map. Let us sketch how this could follow alternatively from Hodge-theoretic considerations. 

Let $\Pi'$ be the fundamental group of $S'$, and $V'_{\mathbb Z}= R_1f'^{an}_\ast \mathbb Z$; one can show that $\Pi \to \Pi'$ is surjective (given that $S$ is Galois over $S'$). 
The section $\sigma$ gives rise to a cocycle in $ H^1(\Pi, V_{\mathbb Z, s} )= H^1(S, V_{\mathbb Z})$, in fact in $H^1(S, V_{\mathbb Z})^{0,0}$. If $\Delta= 0$, it comes from a cocycle $ e\in H^1(\Gamma, V_{\mathbb Z, s})$, whose image $e' \in   H^1(\Pi', V'_{\mathbb Z, s})= H^1(S', V'_{\mathbb Z})$ lies in $H^1(S', V'_{\mathbb Z})^{0,0}\subset H^1(S, V_{\mathbb Z})^{0,0}$. It follows from \eqref{eq12} + \eqref{eq14}  that $\sigma$ comes from a section of $f'$.   
  \end{rems}

\section{A ``Hodge-theoretic proof" of Theorem 1.2.1} 

  Our proof of Theorem 1.2.1 is a generalization of the previous argument for elliptic schemes, relying on the generalization by Mok and To of Shioda's theorem, and on the Hodge-theoretic interpretation of sections (which depends on the Graber-Starr theorem)\footnote{both results go beyond Hodge theory, whence the quotes.}.

  \subsection{General setting} 

\subsubsection{} Let $S$ be a smooth connected algebraic complex variety. For  $s\in S(\mathbb C)$, let us write $\Pi = \pi_1(S^{an}, s)$. Let $$(\ast)\;\;\;0\to V_{\mathbb Z} \to {\bf V}_{\mathbb Z} \to \mathbb Z_S\to 0$$ be an extension of local systems of torsion-free abelian groups on $S^{an}$. We set $$V = V_{\mathbb Z}\otimes \mathbb  Q, \; {\bf V} = {\bf V}_{\mathbb Z}\otimes \mathbb  Q.$$ 
   
   Let $\Gamma$ (\resp $\bf \Gamma$) be the image of $\pi_1(S(\C),s)$ in $GL(V_s)$ (\resp $GL({\bf V}_s)$), and $H$ (\resp $\bf H$) its Zariski-closure. Let $\Delta$ be the kernel of ${\bf \Gamma} \to \Gamma$ and $U$ be the kernel of ${\bf H} \to H$, so that $U$ is a vector group and $\Delta$ is a subgroup of $U(\mathbb Z)$. One has ${\bf\Gamma} = \Delta  \rtimes \Gamma $, and a surjection 
   $ \Pi \surj {\bf\Gamma} $  which induces an injection  
\begin{equation}\label{eq15} H^1({\bf\Gamma}, V_{\mathbb Z, s}) \inj H^1(\Pi, V_{\mathbb Z, s}).\end{equation}
   
   The extension $(\ast)$ is given by a class $e_{\mathbb Z}\in {\rm{Ext}}^1(\mathbb Z_S, V_{\mathbb Z})= H^1(S^{an}, V_{\mathbb Z}) = H^1(\Pi, V_{\mathbb Z, s})$ which comes from $H^1({\bf\Gamma}, V_{\mathbb Z, s})$.
   
   More precisely, let ${\bf v}\in {\bf V}_s$ map to $1\in \mathbb Z$, and let us consider the map $\phi: U(\mathbb Q) \to V_s, \; u\mapsto u({\bf v})-{\bf v}.$ It has the following three properties:
   
    - $\phi$ is independent of the choice of $\bf v$; in particular, one may replace $\bf v$ by $\gamma({\bf v})$, which also maps to $1$ for any $\gamma \in \bf \Gamma$, 
    
    - $\phi$ is {$\bf \Gamma$-equivariant} for the action on $U(\mathbb Q)$ by conjugation: this follows from the previous point,
    
    - $\phi$ is {injective}: the kernel, viewed as a subgroup of $\bf H$, fixes both $\bf v$ and $V_s$, hence is trivial. 
    
  In particular, $V_s \cap \mathbb Q .\phi(\Delta)$ is stable under $\Gamma$. The split exact sequence 
$$0\to \Delta \to {\bf \Gamma} \to \Gamma \to 0,$$ together with the fact that $\Delta$ acts trivially on $V_s$ gives rise to an exact sequence (\cf \eg \cite[\S 12]{DHW} for exactness on the right)
\begin{equation}\label{eq16} 0 \to H^1(\Gamma, V_s) \to H^1({\bf\Gamma}, V_s)\to H^1(\Delta, V_s)^\Gamma \cong Hom_\Gamma(\mathbb Q.\phi(\Delta), V_s)\to 0.\end{equation}

  

\subsubsection{}  This situation occurs in the setting of an abelian scheme $f: A \to S$ and a section $\sigma$ of $f$, with $V_{\mathbb Z} = R_1f^{an}_\ast \mathbb Z  $. 
  In the situation of Theorem  1.2.1, one has to show that the image $e$ of $e_{\mathbb Z}$ in $H^1(\Pi, V_{ s})^{0,0}$ is $0$. On the other hand, by assumption, $e$ comes from a cocycle in  $H^1({\bf \Gamma}, V_{ s})$.

\subsection{Reduction to the case $\Delta =0$} $\,$

 It turns out that one may assume that $A/S$ is simple, but it is not clear a priori that $\mathbb Q . \Delta$ is then either $0$ or $ H_1(A_s, \mathbb Q) $. This is related to the issue of 
{\it non-rigid} abelian schemes, which arises when ${\rm{End}}_S R_1f^{an}_\ast \mathbb Z $ (endomorphisms of the local sytem) is bigger that ${\rm{End}}_S A$ (endomorphism of the variation of Hodge structures) - see \cite{mhS1} for a thorough investigation\footnote{initiated by Deligne and Faltings.} of the issue of non-rigid abelian schemes.
  \subsubsection{} By restricting $S$ one may assume that $A/S$ is polarized. One may also assume that $A/S$ is  simple: indeed, if $A$ is isogeneous to a polarized fibered product $A_1\times_S A_2 \times_S \cdots \times_S A_n$ with simple factors, $V$ becomes a sum $\oplus V_i$ and the problem reduces to the analogous problem for each $V_i$.   
    \subsubsection{} Assume that $A/S$ is simple. In order to show that $\Delta$ is either a lattice in $V_s$ or $0$, it is enough to show that $\mathbb Q.\phi(\Delta)$ is sub-Hodge structure of $V_s$; since $\mathbb Q.\phi(\Delta)$ is $\Gamma$-invariant, it is the fiber at $s$ of a local subsystem $W \subset V$, and this amounts to showing that $W$ is a subvariation of Hodge structures of $V$. 

We use the fact that $(\ast)$ comes from a section $\sigma$ and that $\Gamma$ is a lattice in $H$. The latter implies that there is a weakly special subvariety $S'\subset \mathcal A_{g,n}$ such that the image of the canonical map $S\to  
\mathcal A_{g,n}$ lies in $S'$, and the monodromy of $f' : A' \to S'$ is $\Gamma$ (after `adjustment" of subgroups of finite index). The variation $V$ on $S$ comes from a variation $V'$ on $S'$.  
 One has a cartesian square 
       $$\begin{matrix}    A  &  \stackrel{h}\to  & A'\\   f \downarrow&&  \downarrow {f'} \\ S &  \stackrel{g}\to & S'  \end{matrix}$$ 
which induces a commutative square of $\pi_1$'s
  $$\begin{matrix}    V_s \rtimes \Pi  &  \to  & V_s \rtimes \Gamma \\    \downarrow&&  \downarrow  \\ \Pi &  \to &  \Gamma  \end{matrix}$$ 
 and a commutative diagram with exact colums and injective rows 
   $$\begin{matrix}    H^1(S'^{an}, V')  =  H^1(\Gamma, V_s)  &  \stackrel{g^\ast}{\to} & H^1(  S^{an},   V) = H^1(\Pi, V_s)       \\    {f'}^\ast\downarrow&&  \downarrow f^\ast  \\ H^1(A'^{an}, {f'}^\ast\mathcal V)= H^1(V_s \rtimes \Gamma , V_s) &  \stackrel{h^\ast}{\to} &  H^1(A^{an}, h^\ast V)= H^1( V_s \rtimes\Pi, V_s) \\    \downarrow&&  \downarrow   \\ H^1(V_s, V_s)^\Gamma = {\rm{End}}R_1f_\ast \mathbb Q &\stackrel{\sim}\to &  H^1(V_s, V_s)^\Pi \end{matrix}$$ 
 The section $\sigma$ gives rise to a retraction $H^1( V_s \rtimes\Pi, V_s) \to  H^1(  \Pi, V_s)$ of $f^\ast$, and the composition $\sigma^\ast h^\ast$ factors as
 $$H^1(A'^{an}, {f'}^\ast\mathcal V)= H^1(V_s \rtimes \Gamma , V_s)  \surj H^1({\bf \Gamma}, V_s) \inj  H^1(  S^{an},   V) = H^1(\Pi, V_s).$$ 
 It follows that $H^1({\bf\Gamma} , V_s)$ carries a mixed Hodge structure. The quotient map $$H^1(V_s \rtimes \Gamma, V_s)/H^1({ \Gamma} , V_s)= {\rm{End}}\,R_1f_\ast \mathbb Q  \to H^1({\bf\Gamma} , V_s)/H^1({ \Gamma} , V_s) = Hom_\Gamma(\mathbb Q.\phi(\Delta), V_s)$$ is a morphism of Hodge structures (of weight $0$) and a morphism of ${\rm{End}}\,R_1f_\ast \mathbb Q $-modules. It is thus given by an idempotent element of $({\rm{End}}\,R_1f_\ast \mathbb Q)^{0,0}$, hence $\mathbb Q.\phi(\Delta)$ is a Hodge substructure of $V_s$.
 Since $A/S$ is simple, this is $0$ or $V_s$, hence $\Delta = 0$ or is a lattice in $U_\mathbb R$.  
 


\subsection{Reduction to finiteness of the Mordell-Weil groups of Kuga fiber spaces}  
$\;$

Let us assume that $\Delta = 0$, \ie $\Gamma = \bf\Gamma$. 
  This implies that $ e\in  H^1(S, V)^{0,0}$ comes from  $H^1(S', V')^{0,0}$. 
    By Theorem 3.3.1 (applied to $S'$ instead of $S$), $e$ actually comes from an element of $\mathcal O_{A'}(S')\otimes \mathbb Q$. 
     To conclude, we notice that ${f'} : A' \to S'$ is a Kuga family of polarized abelian varieties in the sense of \cite{MT} (which is more general than the usual definition). The main theorem of \cite{MT}\footnote{which generalizes earlier works by Shioda and Silverberg.} asserts that there is no non-torsion section of $f'$. Hence $e=0$,  $\sigma $ is torsion, and $U=0$.\qed
 
\section{Variants and complements}
\subsection{Vanishing theorems of $H^1$ for arithmetic groups} 
 \subsubsection{The Margulis-Starkov vanishing theorem} Let $H_\R$ be a semisimple algebraic group over $\mathbb R$, and $\Gamma \subset H_\R$ a lattice (for instance, an arithmetic subgroup of a semisimple algebraic group over $\mathbb Q$). Improving on previous results by Margulis, Starkov \cite{St} proved that {\it if there is no epimorphism $\phi : H_\R\to L$ to a Lie group locally isomorphic to $SO(1,n)$ or $SU(1,n)$ such that $\phi(\Gamma)$ is a lattice in $L$, then for any finite-dimensional real representation $W_\mathbb R$ of $\Gamma$, $H^1(\Gamma, W_\mathbb R)=0$}.  
  
\subsubsection{} This applies to the situation $\Gamma \subset H$ of  \S 1.1. In that case, $H_\R$ must be of hermitian type, so that only the case when $L $ is locally isomorphic to $SU(1,n)$ has to be considered (note that the case of \S 2.1 corresponds to $SU(1,1)\cong SL(2)$).  

For instance, the case considered in \cite{DT} is to $H= Sp_{2g}, \, g>1$, so that $H^1(\Gamma, V_s)=0$ by Starkov's theorem, so that our cocycle $e= 0$ and one gets Theorem 1.2.1 as an immediate consequence (without using the Mok-To theorem).

 \subsubsection{Reduction to the case where $H$ is $\mathbb Q$-simple}  One can go a little further in the reduction, assuming, as we may, that $A/S$ is simple and that $H$ is connected. 

 Assume that there is a non-trivial almost direct product decomposition (of semisimple groups) $H= H_1.H_2$. Then there is a subgroup of finite index of $\Gamma$ of the form $\Gamma_1\times \Gamma_2$ where $\Gamma_i$ is a non-trivial arithmetic subgroup of $H_i$. Passing to a finite etale cover of $S$, one may assume that $\Gamma =\Gamma_1\times \Gamma_2$.  
 
Then $V_s^{\Gamma_1} = V_s^{\Gamma_2}= 0$, and $H^1(\Gamma, V_s)= 0$. 

  Indeed, $V_s^{\Gamma_i}$ is $H$-stable; and since $H$ is normal in the Tannakian group $G$ of the variation of Hodge structures $V$, which is connected like $H$ \cite{A1}\cite{A3}, $G$ acts on ${\rm{Aut}}\,H$ by inner automorphisms of $H$, hence normalizes $H_i$. 
Therefore $V_s^{\Gamma_1}, V_s^{\Gamma_2}$ are, up to isogeny, fibers at $s$ of abelian subschemes. Since $\Gamma$ acts faithfully and $A/S$ is simple, one has $V_s^{\Gamma_1} = V_s^{\Gamma_2}= 0$. 
Then the exact sequence $H^1(\Gamma_1, V_s^{\Gamma_2})\to  H^1(\Gamma, V_s)\to  H^1(\Gamma_2, V_s)^{\Gamma_1} = H^1(\Gamma_2, V_s^{\Gamma_1})$ shows that $H^1(\Gamma, V_s)= 0$. 

Therefore one may assume that $H$ is (almost) $\mathbb Q$-simple; it is well-known that it is then a restriction of scalars ${\rm{Res}}_{F/\mathbb Q}\, H^F$, where $H^F$ is an almost absolutely simple group over a totally real field $F$); and further that some real embedding $H_F\otimes_F \mathbb R$ is locally isomorphic to $SU(1,n)$ (the other real embeddings being compact).
 
\subsection{The pro-$\ell$ picture; congruence subgroup property and vanishing theorems for $H^1$}
 \subsubsection{ } Let $\ell$ be a prime number. 
 
 Let $\Gamma_\ell$ (\resp $\bf \Gamma_\ell$) be the image of of the etale fundamental group $\pi_1^{et}(S, s)$ in $GL(V_s\otimes \mathbb Q_\ell)$ (\resp $GL({\bf V}_s \otimes \mathbb Q_\ell)$), an $\ell$-adic analytic Zariski-dense subgroup of $H_{\mathbb Q_\ell}$ (\resp $\bf H_{\mathbb Q_\ell}$). 

Note that {\it derived Lie subalgebras of an algebraic Lie algebra are algebraic} \cite[\S 7, Cor. 7.9]{B}. 
 In particular $ {\rm{Lie}}\,  \Gamma_\ell =  {\rm{Lie}}\,  H_\ell, \; D{\rm{Lie}}\, {\bf\Gamma}_\ell = D{\rm{Lie}}\, \bf H_{\mathbb Q_\ell}$. The same holds for the Lie algebra of an analytic subgroup of a vector group (such as ${\bf H}_{\mathbb Q_\ell}^{ab}$).  One deduces that ${\rm{Lie}}\, {\bf\Gamma}_\ell =  {\rm{Lie}}\, \bf H_\ell$, and further that $\ker ({\bf \Gamma}_\ell  \to \Gamma_\ell)$ is Zariski-dense in $U_{\mathbb Q_\ell}$.
 
  \begin{rem} If $\alpha\in \mathbb Z_\ell^\ast$ is not a root of unity, the  image of $\mathbb Z_\ell  \to GL_{2,\mathbb Q_\ell}, \, z\mapsto   \begin{pmatrix} \alpha & 1\\ 0 &\alpha \end{pmatrix}^z $ is not open in its Zariski-closure $\mathbb G_a \times \mathbb G_m$\footnote{this example also shows that the property of being open in the Zariski-closure cannot be tested on the semisimplification.}.  \end{rem}
 
 \subsubsection{ } Let $\tilde H$ denote the universal covering of $H$, and $\tilde\Gamma$ be the inverse image of $\Gamma$ in $\tilde H(\mathbb Q)$. Note that $\tilde H$ has no compact factor, hence the strong approximation theorem applies to $\tilde H$. If moreover $\tilde H$ satisfies the {\it congruence subgroup property}, Bass-Milnor-Serre show\footnote{their proof uses an argument similar to the one above \S 5.2.1.}  that  $H^1(\tilde\Gamma, W) =0 $ for any finite-dimensional $\mathbb Q$-representation of $\tilde\Gamma$ \cite[Theorem 16.2]{BMS}. In particular, $H^1(\Gamma, V_s)  \subset H^1(\tilde\Gamma, V_s) =0.$

 \subsection{Vista: Integral relative monodromy of normal functions}

\subsubsection{} Let $S$ be a smooth connected complex algebraic variety and  $f: X\to S$ a projective smooth morphism. Let us consider the local system $V_{\mathbb Z} = R^{2n-1}f^{an}_\ast \mathbb Z(n)\,$ (for some $n\in \mathbb N$): it underlies a variation of Hodge structures of weight $-1$, as in 3.2.2.   

Any extension of variations of mixed Hodge structures $0\to V_{\mathbb Z}  \to {\bf V}_{\mathbb Z}  \to \mathbb Z_S\to 0$ corresponds to a holomorphic section (called a horizontal {\it normal function}) of the relative intermediate Jacobian $J^n(X/S) $, which is fibered in complex tori over $S^{an}$: any such extension corresponds to a holomorphic section $\nu$ of $J^n(X/S) $. Standard examples of admissible normal functions come from families of algebraic cycles of codimension $n$ which are homologically trivial in the fibers, by the Abel-Jacobi construction, \cf \eg \cite[\S 2]{Ch}. 

 Theorem 1.2.1 generalizes to this situation whenever the monodromy of $V_{\mathbb Z}$ is a lattice satisfying Starkov's condition: the relative integral monodromy of ${\bf V}_{\mathbb Z} $ is a lattice in the relative algebraic monodromy. 
  
\subsubsection{}  Other examples of admissible normal functions occur in the Green-Griffiths approach to the Hodge conjecture. One starts with a smooth projective variety $\bar X$ of dimension $2n$ and a Hodge cycle $\xi\in H^{2n}(\bar X, \mathbb Z)(n)$. Let $X\to S\subset {\mathbb P}^1$ be constructed from a Lefschetz pencil of hyperplane sections of $\bar X$ (blowing up the axis of the pencil). Assume moreover for simplicity\footnote{otherwise one has to replace  $H^{2n-1}(X_s)(n)$ by the vanishing part of the cohomology} that $H^{2n-1}(\bar X, \mathbb Z)=0$. This provides a normal function $\nu$, \cf \eg \cite[\S 3]{Ch}\footnote{the Green-Griffiths approach relates the algebraicity of $\xi$ to the singularities of $\nu$}. 

In this situation, \cite[\S 6.7]{J} presents mild geometric conditions (conjecturally always satisfied according to \loccit) which guarantee that the monodromy of $V_{\mathbb Z}$ is a lattice in $H_{\mathbb R} = {\rm{Sp}} (V_s \otimes {\mathbb R})$. Theorem 1.2.1 then generalizes to this situation if the lattice is of rank $\geq 4$: the relative integral monodromy of ${\bf V}_{\mathbb Z} $ is of finite index in ${V}_{\mathbb Z} $.

\bigskip\noindent 
\begin{small} \smit{Acknowledgements.} \rm{I am grateful to P. Corvaja and U. Zannier for having brought to my attention the problem of relative integral monodromy versus relative algebraic monodromy, for having explained to me their work in the already problematic case of an elliptic pencil, and for pointing out Example 2.2.1 which prevented a mistake on my side.
 I also thank F. Tropeano for a recent discussion on the topic, which reactivated my interest and led to this note, and G. Baldi for his comments.}  \end{small}

\end{sloppypar}

\begin{thebibliography}{I}
    \bibitem{A1} Andr\'e Y.,  Mumford-Tate groups of mixed Hodge structures and the theorem of the fixed part, Compos. Math.
{82} (1992), 1-24. 
  \bibitem{A2} Andr\'e Y., Groupes de Galois motiviques et p\'eriodes, Exp. 1104 in {\it S\'eminaire Bourbaki 2015/2016}, Ast\'erisque 390 (2017).
    \bibitem{A3} Andr\'e Y., Normality criteria, and monodromy of variations of mixed Hodge structures, appendix to: Kahn B., Albanese kernels and Griffiths groups, Tunisian J. of Math. 3 (2021), 589-656.  
   \bibitem{ACZ} Andr\'e Y., Corvaja P., Zannier U., The Betti map associated to a section of an abelian scheme, Invent. Math. 222 (2020), 161-202.
    \bibitem{Ar} Arapura D.,  The Leray spectral sequence is motivic, math.AG/03010140
   \bibitem{BMS} Bass H., Milnor J., Serre J.-P.,  Solution of the congruence subgroup problem for $SL_n (n \ge 3)$ and $Sp_{2n} (n \ge 2)$. Publ. math. IHES 33 (1967), 59-137. 
 \bibitem{B} Borel A., {\it Linear algebraic groups}, 2nd Ed., G.T.M. 126, Springer-Verlag (1991).
\bibitem{BH} Borel A., Harish-Chandra, Arithmetic subgroups of algebraic groups,
 Annals of Math. 75, 3 (1962), 485-535.
 \bibitem{BFZP} P. Brosnan, Fang H., Nie Z., Pearlstein G., Singularities of admissible normal functions, Invent. Math  177 (2009), 599-629.
 \bibitem{Ca} Carlson J., The geometry of the extension class of a mixed Hodge structure, in {\it Algebraic Geometry}, Bowdoin, 1985, Proc. Symposia in Pure Math. 46, Amer. Math. Soc. (1987), 199-222.
 \bibitem{CK} Cattani E., Kaplan A., Degenerating variations of Hodge structure, in {\it{Th\'eorie de Hodge}} (Luminy), Ast\'erisque, 179-180.
   \bibitem{Ch} Charles F., Progr\`es r\'ecents sur les fonctions normales d'apr\`es Green-Griffiths, Brosnan-Pearlstein, M. Saito, Schnell..., Exp. 1063 in {\it S\'eminaire Bourbaki 2014}, Ast\'erisque 361 (2014).
 \bibitem{CZ} Corvaja P., Zannier U., Unramified sections of the Legendre scheme and modular forms, J.
Geom. Phys., 166, 26 (2021). Corrigendum (2024). 
    \bibitem{DHW} Dekimpe K., Hartl M., Wauters S., A seven-term exact sequence for the cohomology of group extension, J. Algebra 369 (2012).
      \bibitem{D1} Deligne P., Th\'eorie de Hodge II, Publ. Math. I.H.E.S. 40 (1971), 5-57.
   \bibitem{D} Deligne P., Th\'eorie de Hodge III, Publ. Math. I.H.E.S. 44 (1974), 5-77.
 \bibitem{DT} Dolce P., Tropeano F., Relative monodromy of abelian logarithms for finite covers of universal families, http://arxiv.org/pdf/2407.19476
  \bibitem{DT2} Dolce P., Tropeano F., Relative monodromy of abelian logarithms of ramified sections on abelian schemes, preprint (2024).
       \bibitem{F} Fakhruddin N.,  Restriction of sections of abelian schemes, arXiv:math/0310405v1
     \bibitem{GS} Graber T., Starr J., Restriction of Sections for Families of Abelian
Varieties, in {\it A Celebration of Algebraic Geometry: A Conference in Honor of
J. Harris}, Clay Math. Proc., volume 18.
      \bibitem{J} Janssen W., Skew-symmetric vanishing lattices and their monodromy groups, Math. Ann. 266 (1983), 115-133.
       \bibitem{K} Kashiwara M.,  A study of variation of mixed Hodge structure, Publ. RIMS 22 (1986), 991-1024.
              \bibitem{KP} Kerr M., Pearlstein G, An exponential history of functions with logarithmic growth, in {:it Topology of Stratified Spaces}
MSRI Publications
 58 (2011).
         \bibitem{MT} Mok N., To W.K., Eigensections on Kuga families of abelian varieties and finiteness of their Mordell-Weil groups, J. reine und angew. Mathematik, 444 (1993), 29–78.
     \bibitem{N} Nori M., A nonarithmetic monodromy group. C. R. Acad. Sci. Paris  Math. 302(2), 71-72 (1986).
       \bibitem{mS2} Saito M., Mixed Hodge modules, Publ. R. I. M. S. 26 (1990), 221-333.
     \bibitem{mS2} Saito M., Admissible normal functions, J. Alg. Geom. 5 (1996), 235-276.
     \bibitem{mhS1} Saito M.-H., Classification of non-rigid families of abelian varieties, Tohoku Math. J. 45 (1993), 159-189.
       \bibitem{mhS2} Saito M.-H., Finiteness of Mordell-Weil groups of Kuga fiber spaces of abelian varieties. Publ. RIMS, 29 (1993), 29-62.
     \bibitem{Sh}  Shioda T., On elliptic modular surfaces, Journal of the M.S.J, 24 (1972), 20-59.
      \bibitem{St} Starkov A., Vanishing of the first cohomology for lattices in Lie groups, Journal of Lie Theory
Volume 12 (2002), 449-460.
          \bibitem{SZ} Steenbrink J., Zucker S.,  Variation of mixed Hodge structure, I, Invent. Math.
80 (1985), 489-542.
    \bibitem{V}  Voisin C., {Torsion points of sections of Lagrangian torus fibrations
and the Chow ring of hyper-K\"ahler manifolds},  in {\it Geometry of Moduli}, Abel Symposia, Springer (2018), 295-326.
  \bibitem{Z} Zucker S, Hodge theory with degenerating coefficients, Ann. of Math.
 109 (1979), 415-476.
 
   
    \end{thebibliography}
\end{document}